\documentclass{cmslatex}
\usepackage{latexsym, amssymb, enumerate, amsmath}

\usepackage[dvips]{graphicx}

\renewcommand{\theequation}{\arabic{section}.\arabic{equation}}

\def\half{\frac 1 2}

\renewcommand{\epsilon}{\varepsilon}

\newcommand{\bc}{\mathbb{C}}
\newcommand{\br}{\mathbb{R}}

\newcommand{\bA}{\boldsymbol{A}}
\newcommand{\balpha}{\boldsymbol{\alpha}}
\newcommand{\bphi}{\boldsymbol{\Phi}}
\newcommand{\bpsi}{\boldsymbol{\Psi}}
\newcommand{\bu}{\boldsymbol{u}}
\newcommand{\bvarphi}{\boldsymbol{\varphi}}
\newcommand{\bx}{\boldsymbol{x}}
\newcommand{\bxi}{\boldsymbol{\xi}}
\newcommand{\by}{\boldsymbol{y}}

\newcommand{\rd}{\mathrm{d}}
\newcommand{\mca}{\mathcal{A}}
\newcommand{\mcd}{\mathcal{D}}
\newcommand{\mcl}{\mathcal{L}}
\newcommand{\nn}{\nonumber}
\newcommand{\nin}{\noindent}

\newcommand{\abs}[1]{\left\vert#1\right\vert}
\newcommand{\brac}[1]{\left(#1\right)}
\newcommand{\norm}[1]{\left\Vert#1\right\Vert}

\newtheorem{thm}{Theorem}[section]

\newtheorem{rem}[thm]{Remark}
\begin{document}

\title{Gaussian Beam Methods for the Dirac Equation in the Semi-classical Regime
\thanks{{Received date / Revised version date.}
          H. Wu was supported by the NSFC project 11101236. Z.Y. Huang was support by the NSFC project 11071139 and the National Basic Research Program of China under the grant 2011CB309705. S. Jin was support by the NSF grants DMS-0608720, DMS-1114546 and NSF FRG grant DMS-0757285. S. Jin was also supported by a Van Vleck Distinguished Research Prize and a Vilas Associate Award from the University of Wisconsin-Madison. D.S. Yin was supported by the NSFC project 10901091.
}}
          %For each author, make a block with the following four macros:
\author{Hao Wu\thanks {Department of Mathematical Sciences, Tsinghua University, Beijing, 100084, China, (hwu@tsinghua.edu.cn).}
\and Zhongyi Huang\thanks {Department of Mathematical Sciences, Tsinghua University, Beijing, 100084, China, (zhuang@math.tsinghua.edu.cn).} \and
Shi Jin\thanks{Department of Mathematics, Institute of Natural Sciences and Key Lab of
Scientific and Engineering Computing-Ministry of Education, Shanghai Jiao Tong University, Shanghai 200240, P.R. China; and Department of Mathematics, University of Wisconsin-Madison, Madison, WI 53706, USA, (jin@math.wisc.edu).}
\and Dongsheng Yin\thanks {Department of Mathematical Sciences, Tsinghua University, Beijing, 100084, China, (dyin@math.tsinghua.edu.cn).}}
          %{Put the URL for your home page here if you have one}

          %Use \thanks statements for acknowledgements of grants and
          %support. They will appear below all the authors' addresses, so be
          %specific about which author is thanking whom:

          %\thanks{}
          
\pagestyle{myheadings} \markboth{Gaussian beam for the Dirac equation}{Hao Wu, Zhongyi Huang, Shi Jin and Dongsheng Yin}\maketitle

\begin{abstract}
The Dirac equation is an important model in relativistic quantum mechanics. In the semi-classical regime $\epsilon\ll1$, even a spatially spectrally accurate time splitting method \cite{HuJi:05} requires the mesh size to be $O(\epsilon)$, which makes the direct simulation extremely expensive. In this paper, we present the Gaussian beam method for the Dirac equation. With the help of an eigenvalue decomposition, the Gaussian beams can be independently evolved along each eigenspace and summed to construct an approximate solution of the Dirac equation. Moreover, the proposed Eulerian Gaussian beam keeps the advantages of constructing the Hessian matrices by simply using level set functions' derivatives. Finally, several numerical examples show the efficiency and accuracy of the method.
\end{abstract}


\section{Introduction}\label{sec:intro}
We are interested in developing the Gaussian beam method for the Dirac equation in the semi-classical regime
\begin{equation} \label{eqn:dirac}
    i\epsilon\partial_t \bpsi^{\epsilon}=-i\epsilon\balpha\cdot \nabla\bpsi^{\epsilon}-
        \balpha\cdot \bA\bpsi^{\epsilon}+\beta\bpsi^{\epsilon}+V\bpsi^{\epsilon},
\end{equation}
subject to the Cauchy initial data
\begin{equation} \label{eqn:diracinit}
    \bpsi^{\epsilon}(0,\bx)=\bu_I(\bx)^{iS_I(\bx)/\epsilon}, \quad \bx\in\br^3.
\end{equation}
Here $\bpsi^{\epsilon}(t,x)=(\Psi^{\epsilon}_1, \Psi^{\epsilon}_2, \Psi^{\epsilon}_3, \Psi^{\epsilon}_4)^T\in \bc^4$ is the spinor field, normalized s.t.,
\begin{equation*}
    \int_{\br^3}\abs{\bpsi^{\epsilon}(t,\bx)}^2 \rd \bx=1,
\end{equation*}
$0<\epsilon\ll1$ denotes the semi-classical parameter, $V(t,\bx)\in\br$ is the external electric potential and $\bA(t,\bx)\in\br^3$ represents the external magnetic potential, i.e., $\bA=(A_1,A_2,A_3)$. Without loss of generality, we only consider static external field in this paper. The Dirac matrices $\beta,\; \balpha=(\alpha_1,\alpha_2,\alpha_3)$ are complex-valued Hermitian matrices, which are given by
\begin{equation*}
    \beta := \brac{\begin{array}{cc} I_2 & 0 \\ 0 & -I_2\end{array}}, \quad
    \alpha^k :=\brac{\begin{array}{cc} 0 & \sigma^k \\ \sigma^k & 0 \end{array}},
\end{equation*}
with $I_2$ the $2\times2$ identity matrix and $\sigma^k$ the $2\times2$ Pauli matrices, i.e.,
\begin{equation*}
    \sigma^1 :=\brac{\begin{array}{cc} 0 & 1 \\ 1 & 0 \end{array}}, \quad
    \sigma^2 :=\brac{\begin{array}{cc} 0 & -i \\ i & 0 \end{array}}, \quad
    \sigma^3 :=\brac{\begin{array}{cc} 1 & 0 \\ 0 & -1 \end{array}}.
\end{equation*}
The physical observables can be defined in terms of $\; \bpsi^{\epsilon}(t,\bx)$:
\begin{eqnarray}
    \textrm{Particle-density} && \quad \rho^{\epsilon}=\abs{\bpsi^{\epsilon}}^2, \\
    \textrm{Current-density} && \quad j_k^{\epsilon}=\langle\bpsi^{\epsilon},\alpha^k\bpsi^{\epsilon}\rangle_{\bc^4}.
\end{eqnarray}

The Dirac equation \cite{Di:28, Th:92} is a relativistic wave equation which plays a fundamental role in relativistic quantum mechanics. It provides a natural description of the  particles with spin $1/2$ , i.e. electrons, neutrinos, muons, protons, neutrons, etc. The Dirac equation also predicts some peculiar effects, such as Klein's paradox \cite{Kl:29} and ``Zitterbewegung" an unexpected quivering motion of a free relativistic quantum particle \cite{Sc:30}. Recently, graphene \cite{CsDa:06, NeGu:09} and topological insulators \cite{HsQiWr:08, XiQi:09} are studied widely in connection to the Dirac equation in the semiclassical regime \cite{MoSc:11}.

In the semi-classical regime $\epsilon\ll1$, the solution to the Dirac equation is highly oscillatory. Thus, for any domain-based discretization method, the number of mesh points in each spatial direction should be at least $O(\epsilon^{-1})$ \cite{JoBlSa:88}. If the potential is sufficiently smooth, and the initial data of the Dirac equation is compactly supported, the time-splitting spectral method \cite{HuJi:05} offers the best numerical resolution. The spatial meshing strategy is almost of optimal order $O(\epsilon^{-1})$ and the time step can be $O(1)$.

One alternative efficient numerical approach for solving the Dirac equation is the WKB method \cite{SpMa03:44, SpMa03:45, Sp:00}. For a first order approximation, this method tries to seek an asymptotic solution:
\begin{equation} \label{eqn:asy}
     \bpsi^{\epsilon}(t,\bx)=\bu(t,\bx)e^{iS(t,\bx)/\epsilon}+O(\epsilon), \ S\in \br,
\end{equation}
where the amplitude $u$ and the phase $S$ are smooth functions independent of $\epsilon$. Substituting \eqref{eqn:asy} into the Dirac equation, one derives the eikonal equation and the transport equation. Since the eikonal equation is of the Hamilton-Jacobi type, the solution becomes singular after caustic formulation. Beyond caustics, the correct semi-classical solution of the Dirac equation contains several phases. In the last decades, many approaches have been proposed to capture this multi-phased solutions, see reviews \cite{EnRu:03, JiMaSp:11}.

A serious drawback of the WKB method is that the solution ceases to be valid at caustics where the rays intersect and the amplitudes blow up. The Gaussian beam method, which was first proposed by Heller in quantum chemistry \cite{He:75} and independently developed by Popov in Geophysics \cite{Po:82}, is an efficient approach that allows accurate computation of the amplitude and phase information near caustics. The main difference between the WKB method and the Gaussian beam method is that the Gaussian beam method allows the phase function to be complex off the center of the beam and the imaginary part of the phase function is positive definite, which makes the solution decay exponentially away from the center. The validity of the Gaussian beam method at caustics was analyzed by Ralston in \cite{Ra:82}. The Gaussian beam and related methods have become very popular in high frequency waves problems \cite{LiRuTa:pre, LuYa:11, LuYa:CPAM, MoRu:10, QiYi:10, TaEnTs:09, Ta:08, YiZh:11, YiZh:12} in recent years. Most of the methods were in the Lagrangian framework. More recently, Eulerian Gaussian beam methods have also received a special attention for its advantages of uniform accuracy \cite{JiWuYa:08, JiWuYa:10, JiWuYa:11, JiWuYaHu:10, LeQi:09, LeQiBu:07}. A major simplification of the Eulerian Gaussian beam method is that the Hessian matrices can be constructed by taking derivatives of the level set functions \cite{JiWuYa:08}. This greatly reduces the computational cost.

To our knowledge, no Gaussian beam methods have previously been developed for the Dirac equation. It is the goal of this paper to develop such a method by extending the previous method of \cite{JiWuYa:08} for the Schr\"odinger equation to the Dirac equation \eqref{eqn:dirac}-\eqref{eqn:diracinit}. With the help of the eigenvalue decomposition, the Gaussian beams evolve independently of each other. Moreover, the energy transition is forbidden since the Dirac matrix $\beta$ results an $O(1)$ band gap. Being different from the Gaussian beam methods for the Schr\"odinger equation, the higher order Taylor expansion and asymptotic expansion must be considered for the amplitude. After making use of the solvability condition and matching the different expansions, one gets the transport equation for the lower order term of the amplitude. When the evolution is done, the solution can be simply constructed by the summation of all Gaussian beams. The solution will be shown to have a good accuracy even around caustics, with a coarse mesh size of $O(\sqrt{\epsilon})$ and large time step of $O(\sqrt{\epsilon})$. A remarkable aspect of the Eulerian Gaussian beam method is that it still possesses the advantage of the previous method \cite{JiWuYa:08}, which is an important benefit for the $3D$ simulation.

The paper is organized as follows. After reviewing the semi-classical limit of the Dirac equation in Section \ref{sec:DIC}, we formulate the Gaussian beam method for \eqref{eqn:dirac}-\eqref{eqn:diracinit} in Section \ref{sec:GBM}. In Section \ref{sec:numexa}, our method is shown to be accurate and efficient by several numerical examples. Finally, we conclude the paper in Section \ref{sec:con}.
\medskip

%*****************************************************************************************
%************* Dirac *********************************************************************
%*****************************************************************************************

\section{The Dirac equation and the semi-classical limit}\label{sec:DIC}
Denote the Dirac operator by
\begin{equation*}
    \mcd\brac{\bx,\bxi} = \balpha\cdot \brac{\bxi-\bA(\bx)}+\beta+V(\bx),
\end{equation*}
then we have
\begin{equation*}
    \mcd\brac{\bx,-i\epsilon\nabla}\bpsi^{\epsilon}=
        \balpha \cdot\brac{-i\epsilon\nabla-\bA(\bx)}\bpsi^{\epsilon}+\beta\bpsi^{\epsilon}+V(\bx)\bpsi^{\epsilon}.
\end{equation*}
Therefore the semi-classically scaled Dirac equation \eqref{eqn:dirac} can be written as
\begin{equation*}
    i\epsilon\partial_t \bpsi^{\epsilon}=\mcd\brac{\bx,-i\epsilon\nabla}\bpsi^{\epsilon}.
\end{equation*}

Let
\begin{equation*}
    \lambda(\bx,\bxi)=\sqrt{\abs{\bxi-\bA(\bx)}^2+1},
\end{equation*}
then
\begin{equation*}
    h^{\pm}(\bx,\bxi) = \pm \lambda(\bx,\bxi)+V(\bx)
\end{equation*}
are two different eigenvalues, each with multiplicity two, of the Dirac operator $\mcd\brac{\bx,\bxi}$. The corresponding projectors $\Pi^{\pm}(\bx,\bxi)$ are given by
\begin{equation*}
    \Pi^{\pm}(\bx,\bxi)=\half\brac{I_4 \pm \frac{1}{\lambda(\bx,\bxi)}
        \brac{\mcd\brac{\bx,\bxi}-V(\bx)I_4}}.
\end{equation*}

Plugging the following WKB-ansatz into \eqref{eqn:dirac},
\begin{equation*}
    \bpsi^{\epsilon}(t,\bx)=e^{iS(t,\bx)/\epsilon}\sum_{j=0}^{\infty}\epsilon^j\bu_j(t,\bx),
\end{equation*}
where $\bu_j(t,\bx)\in C^{\infty}(\br^4,\bc^4)$, matching the $O(1)$ and $O(\epsilon)$ asymptotic coefficients, one has
\begin{eqnarray}
    \brac{\partial_t S+\mcd(x,\nabla S)}\bu_0 &=& 0, \label{eqn:asy0} \\
    i\brac{\partial_t+\balpha\cdot \nabla}\bu_0-\brac{\partial_t S+\mcd(\bx,\nabla S)}\bu_1 &=& 0. \label{eqn:asy1}
\end{eqnarray}
In order to get a nontrivial solution $\bu_0(t,\bx)\ne0$ in \eqref{eqn:asy0}, one gets
\begin{equation*}
    \textrm{det}\brac{\partial_tS+\mcd(\bx,\nabla S)}=0,
\end{equation*}
which leads to the eikonal equation
\begin{equation} \label{eqn:eikonal}
    \partial_t S^{\pm}+h^{\pm}\brac{\bx,\nabla S^{\pm}}=0.
\end{equation}
Applying the projection $\Pi^{\pm}(\bx,\nabla S^{\pm})$ to \eqref{eqn:asy1}, one gets the solvability condition
\begin{equation*}
    \Pi^{\pm}(\bx,\nabla S^{\pm})\brac{\partial_t+\balpha\cdot \nabla}\bu_0^{\pm}=0.
\end{equation*}
After a series of calculations \cite{SpMa03:44, SpMa03:45}, the following transport equation can be derived:
\begin{equation} \label{eqn:transport}
    \partial_t\bu_0^{\pm}+\brac{\omega^{\pm}(\bx,\nabla S^{\pm})\cdot\nabla}\bu_0^{\pm}
        +\half\brac{\nabla\cdot\omega^{\pm}(\bx,\nabla S^{\pm})}\bu_0^{\pm}=\mca^{\pm}(\bx,\nabla S^{\pm})\bu_0^{\pm},
\end{equation}
with
\begin{eqnarray*}
    \omega^{\pm}(\bx,\bxi) &=& \nabla_{\bxi}h^{\pm}(\bx,\bxi), \\
    \mca^{\pm}(\bx,\bxi) &=& \sum_{k\ne l}\frac{\alpha^k\alpha^l}{2\lambda}\brac{\partial_{x_k}A_l}
        -\frac{1}{2\lambda}\balpha\cdot\nabla h^{\pm}-
        \frac{1}{2\lambda^2}\balpha\cdot\brac{\brac{\bxi-\bA}\cdot \nabla\bA} \nn \\
    && +\frac{1}{2\lambda^3}\brac{\bxi-\bA}\cdot\brac{(\bxi-\bA)\cdot\nabla\bA+\lambda\nabla h^{\pm}}.
\end{eqnarray*}
\begin{rem}
    If the external magnetic potential is zero, i.e. $\bA=0$, then
    \begin{equation*}
        \mca^{\pm}(\bx,\bxi)= -\frac{1}{2\lambda}\balpha\cdot\nabla V(\bx)+\frac{1}{2\lambda^2}\bxi\cdot\nabla V(\bx).
    \end{equation*}
\end{rem}
\medskip

%*****************************************************************************************
%************* Gaussian beam *************************************************************
%*****************************************************************************************

\section{The Gaussian beam method}\label{sec:GBM}\setcounter{equation}{0}
In this section, we derive the Gaussian beam method using both Lagrangian and Eulerian formulations. We first introduce the Lagrangian Gaussian beam method for solving the Dirac equation, then discuss the Eulerian Gaussian beam method.

\subsection{The Lagrangian formulation}
In this subsection, we describe how to solve the Dirac equation \eqref{eqn:dirac} by the Lagrangian Gaussian beam method, which is given by the following ansatz:
\begin{equation} \label{eqn:LGB}
    \boldsymbol{\phi}^{\epsilon\pm}_{la}(t,\bx,\by_0)=\bu_0^{\pm}(t,\by^{\pm})e^{iT^{\pm}(t,\bx,\by^{\pm})/\epsilon},
\end{equation}
where $\by^{\pm}=\by^{\pm}(t,\by_0)$ and $T^{\pm}(t,\bx,\by^{\pm})$ is a second order Taylor truncated phase function
\begin{equation*}
    T^{\pm}(t,\bx,\by^{\pm})=S^{\pm}(t,\by^{\pm})+\bxi^{\pm}(t,\by^{\pm})\cdot(\bx-\by^{\pm})
        +\half (\bx-\by^{\pm})^TM^{\pm}(t,\by^{\pm})(\bx-\by^{\pm}).
\end{equation*}
Here $S^{\pm}\in\br,\;\bxi^{\pm}\in\br^3,\;M^{\pm}\in\bc^{3\times3}$ and $\bu_0^{\pm}\in\bc^4$. Then the evolutionary ODEs of the Lagrangian Gaussian beam \eqref{eqn:lgby}-\eqref{eqn:lgbm} can be derived by the eigenvalue decomposition technique and the standard Gaussian beam method (for details see Appendix):
\begin{eqnarray}
    \frac{\rd \by^{\pm}}{\rd t} &=& \nabla_{\bxi}h^{\pm}, \label{eqn:lgby} \\
    \frac{\rd \bxi^{\pm}}{\rd t} &=& -\nabla_{\by}h^{\pm}, \label{eqn:lgbxi}\\
    \frac{\rd  S^{\pm}}{\rd t} &=& \nabla_{\bxi}h^{\pm}\cdot\bxi^{\pm}-h^{\pm}, \label{eqn:lgbs} \\
    \frac{\rd M^{\pm}}{\rd t} &=& -\nabla_{\by\by}h^{\pm}-\nabla_{\by\bxi}h^{\pm}M^{\pm}
        -M^{\pm}\nabla_{\bxi\by}h^{\pm}-M^{\pm}\nabla_{\bxi\bxi}h^{\pm}M^{\pm}, \label{eqn:lgbm}\\
    \frac{\rd \bu_0^{\pm}}{\rd t} &=& -\half \brac{\nabla_{\by}\cdot\omega^{\pm}}\bu_0^{\pm}+
        \mca \bu_0^{\pm}, \label{eqn:lgbu}
\end{eqnarray}
in which $\by^{\pm}=\by^{\pm}(t,\by_0),\;\bxi^{\pm}=\bxi^{\pm}(t,\by^{\pm}(t,\by_0)),\; S^{\pm}=S^{\pm}(t,\by^{\pm}(t,\by_0)),\; M^{\pm}=M^{\pm}(t,\by^{\pm}(t,\by_0))$ and $\bu_0^{\pm}=\bu_0^{\pm}(t,\by^{\pm}(t,\by_0))$. The ODE of Lagrangian Gaussian beam amplitude \eqref{eqn:lgbu} can be given by the solvability condition and high order Gaussian beam formulation. Equations \eqref{eqn:lgby}-\eqref{eqn:lgbxi} are the ray tracing equations, and the equation \eqref{eqn:lgbm} is a Riccati equation which can be alternatively solved by the following dynamic first order system
\begin{eqnarray}
    \frac{\rd P^{\pm}}{\rd t} &=& \brac{\nabla_{\bxi\by}h^{\pm}}P^{\pm}+\brac{\nabla_{\bxi\bxi}h^{\pm}}R^{\pm}, \label{eqn:lgbp} \\
    \frac{\rd R^{\pm}}{\rd t} &=& -\brac{\nabla_{\by\by}h^{\pm}}P^{\pm}-\brac{\nabla_{\by\bxi}h^{\pm}}R^{\pm}. \label{eqn:lgbr}
\end{eqnarray}
Then the Hessian matrices satisfy $M^{\pm}=R^{\pm}\brac{P^{\pm}}^{-1}$. After that, the Lagrangian Gaussian beam solution to the Dirac equation \eqref{eqn:dirac} is constructed as
\begin{equation} \label{eqn:LGB_int}
    \bphi^{\epsilon}_{la}(t,\bx)=\brac{\frac{1}{2\pi\epsilon}}^{\frac{3}{2}}\int_{\br^3}
        \brac{r_{\theta}(\bx-\by^+)\boldsymbol{\phi}^{\epsilon+}_{la}(t,\bx,\by_0)
        +r_{\theta}(\bx-\by^-)\boldsymbol{\phi}^{\epsilon-}_{la}(t,\bx,\by_0)}\rd \by_0,
\end{equation}
where $r_{\theta}\in C_0^{\infty}\brac{\br^3},\;r_{\theta}\ge0$ is a truncation function with $r_{\theta}\equiv1$ in a ball of radius $\theta>0$ about the origin. The discrete form of \eqref{eqn:LGB_int} is given as
\begin{multline} \label{eqn:LGB_dis}
    \bphi^{\epsilon}_{la}(t,\bx)=
        \brac{\frac{1}{2\pi\epsilon}}^{\frac{3}{2}}\sum_{j=1}^{N_{\by_0}}
        r_{\theta}(\bx-\by^+(t,\by_0^j))\boldsymbol{\phi}^{\epsilon+}_{la}(t,\bx,\by_0^j)
        \Delta \by_0 \\
        +\brac{\frac{1}{2\pi\epsilon}}^{\frac{3}{2}}\sum_{j=1}^{N_{\by_0}}
        r_{\theta}(\bx-\by^+(t,\by_0^j))\boldsymbol{\phi}^{\epsilon+}_{la}(t,\bx,\by_0^j)
        \Delta \by_0,
\end{multline}
where $\by_0^j$ are the equidistant mesh points, and $N_{\by_0}$ is the number of the beams initially centered at $\by_0^j$. The initial conditions for the equations \eqref{eqn:lgby}-\eqref{eqn:lgbu} are as follows \cite{JiWuYa:08, Ta:08}
\begin{eqnarray}
    \by^{\pm}(0,\by_0) &=& \by_0, \label{eqn:lgbyinit}\\
    \bxi^{\pm}(0,\by_0) &=& \nabla S_I(\by_0), \label{eqn:lgbxiinit}\\
    S^{\pm}(0,\by_0) &=& S_I(\by_0), \label{eqn:lgbsinit}\\
    M^{\pm}(0,\by_0) &=& \nabla^2 S_I(\by_0)+iI, \label{eqn:lgbminit}\\
    \bu_0^{\pm}(0,\by_0) &=& \Pi^{\pm}(\by_0,\nabla S_I(\by_0))\bu_I(\by_0). \label{eqn:lgbuinit}
\end{eqnarray}

\begin{rem}
    The evolutionary equation \eqref{eqn:lgbu} and the initial condition \eqref{eqn:lgbuinit} ensure that
    \begin{align*}
        & \Pi^+\bu_0^+(t,\by_0)=\bu_0^+(t,\by_0), \quad \Pi^+\bu_0^-(t,\by_0)=0, \\
        & \Pi^+\bu_0^-(t,\by_0)=0, \quad \Pi^-\bu_0^-(t,\by_0)=\bu_0^-(t,\by_0),
    \end{align*}
    for $\forall t\ge0$. The related discussion for the semi-classical limit can be found in \cite{SpMa03:44}.
\end{rem} \medskip

\begin{rem}
    As discussed in \cite{JiWuYa:08}, to compute the Gaussian beam solutions for the Schr\"odinger equation, the optimal mesh can be $O(\epsilon^{\half})$ and the time step requirement is of order $O(\epsilon^{\half})$. Therefore, the total computational cost is $O(\epsilon^{-2})$ for the $3$D simulation. On the other hand, the computational cost for the direct numerical methods should be at least $O(\epsilon^{-3})$ for accurate physical observables and $O(\epsilon^{-4})$ for accurate wave fields. It is obvious to see that the Gaussian beam method is much more efficient in the semi-classical regime $\epsilon\ll1$. Since the solutions of the Dirac equation has the similar high frequency structures to the ones of the Schr\"odinger equation. The similar discussions can be made to the Dirac equation.
\end{rem}

\subsection{The Eulerian formulation}
In this subsection, the Eulerian Gaussian beam method using the level set method \cite{JiLiOsTs:05, JiOs:03} is introduced to solve the Dirac equation \eqref{eqn:dirac}. The analogous derivations are given in details by the former work of Jin {\it et al} \cite{JiWuYa:08}.

Define the Liouville operator as
\begin{equation*}
    \mcl^{\pm}=\partial_t+\nabla_{\bxi}h^{\pm}\cdot \nabla_{\by}-\nabla_{\by}h^{\pm}\cdot\nabla_{\bxi},
\end{equation*}
then the level set equations corresponding to equations \eqref{eqn:lgby}-\eqref{eqn:lgbu} are given by
\begin{eqnarray}
    \mcl^{\pm}\bvarphi^{\pm} &=& 0, \label{eqn:egbphi} \\
    \mcl^{\pm}S^{\pm} &=& \nabla_{\bxi}h^{\pm}\cdot\bxi^{\pm}-h^{\pm}, \label{eqn:egbs} \\
    \mcl^{\pm}\bu_0^{\pm} &=& -\half (\nabla_{\by}\cdot\omega^{\pm})\bu_0^{\pm}+\mca \bu_0^{\pm}. \label{eqn:egbu}
\end{eqnarray}
Here $\bvarphi^{\pm}=\bvarphi^{\pm}(t,\by,\bxi)\in \bc^3$ are the level set functions. The zero level set of $\textrm{Re}[\varphi_k^{\pm}]$ gives the (multi-valued) velocity. We also have the phase $S^{\pm}=S^{\pm}(t,\by,\bxi)\in\br$ and the amplitude $\bu_0^{\pm}=\bu_0^{\pm}(t,\by,\bxi)\in\bc^4$ in the phase space. To be compatible with the initial data \eqref{eqn:lgbyinit}-\eqref{eqn:lgbuinit}, we use the following initial condition:
\begin{eqnarray}
    \bvarphi^{\pm}(0,\by,\bxi) &=& -i\by+(\bxi-\nabla_{\by}S_i(\by)), \label{eqn:egbphiinit} \\
    S^{\pm}(0,\by,\bxi) &=& S_I(\by), \label{eqn:egbsinit} \\
    \bu_0^{\pm}(0,\by,\bxi) &=& \Pi^{\pm}(\by_0,\nabla S_I(\by_0))\bu_I(\by_0). \label{eqn:egbuinit}
\end{eqnarray}
From \eqref{eqn:egbphi} and \eqref{eqn:egbphiinit}, the Hessian matrices are constructed via
\begin{equation*}
    M^{\pm}=-\nabla_y\bvarphi\brac{\nabla_{\bxi}\bvarphi}^{-1}.
\end{equation*}
As a result of this property, we don't need to solve the level set equations for $M^{\pm},\;P^{\pm}$ or $R^{\pm}$ corresponding to equations \eqref{eqn:lgbm} and \eqref{eqn:lgbp}-\eqref{eqn:lgbr} as was done in the Eulerian Gaussian beam method \cite{LeQi:09, LeQiBu:07} . This can save a lot of computational resources especially for such $3$D problem.
After that, the Eulerian Gaussian beam solution to the Dirac equation \eqref{eqn:dirac} is constructed as
\begin{multline} \label{eqn:EGB_int}
    \bphi^{\epsilon}_{eu}(t,\bx)=\brac{\frac{1}{2\pi\epsilon}}^{\frac{3}{2}}\int_{\br^3}\int_{\br^3}
        r_{\theta}(\bx-\by)\Big(\boldsymbol{\phi}^{\epsilon+}_{eu}(t,\bx,\by,\bxi)
        \delta(\textrm{Re}[\bvarphi^+(t,\by,\bxi)]) \\
        +\boldsymbol{\phi}^{\epsilon-}_{eu}(t,\bx,\by,\bxi)
        \delta(\textrm{Re}[\bvarphi^-(t,\by,\bxi)])\Big)
        \rd \bxi \rd \by,
\end{multline}
where
\begin{eqnarray*}
    \boldsymbol{\phi}^{\epsilon\pm}_{eu}(t,\bx,\by,\bxi) &=&
        \bu_0^{\pm}(t,\by,\bxi)e^{iT^{\pm}(t,\bx,\by,\bxi)/\epsilon}, \\
    T^{\pm}(t,\bx,\by,\bxi) &=& S^{\pm}(t,\by,\bxi)+\bxi\cdot(\bx-\by)
        +\half (\bx-\by)^TM^{\pm}(t,\by,\bxi)(\bx-\by).
\end{eqnarray*}

\begin{rem}
    The equation \eqref{eqn:EGB_int} can be solved by a discretized delta function integral method \cite{Wen:09, Wen:10} or a local semi-Lagrangian method \cite{JiWuYa:08, LeQiBu:07}.
\end{rem}
\medskip

\begin{rem}
    One can use the local level set method \cite{Mi:03, PeMeOs:99} to solve the level set equation in the vicinity of a lower-dimensional zero level curve of $\textrm{Re}[\varphi^{\pm}]$ to reduce the total computational cost for the Eulerian Gaussian beam method. An alternative efficient way is to use the semi-Eulerian Gaussian beam method proposed in \cite{JiWuYa:11}.
\end{rem}
\medskip
%*****************************************************************************************
%************* Numerical examples ********************************************************
%*****************************************************************************************

\section{Numerical examples}\label{sec:numexa}\setcounter{equation}{0}
In this section, we present several numerical examples to show the accuracy and efficiency of the Gaussian beam method. We compute the solution of the Dirac equation \eqref{eqn:dirac} by the time-splitting spectral scheme \cite{HuJi:05}. The reference solutions $\bpsi^{\epsilon}(t,\bx)$ are computed on a very fine mesh and a very small time step. In all the numerical examples, the truncation parameter $\theta$ in \eqref{eqn:LGB_dis} is chosen large enough so that the cut-off error is almost zero. \bigskip

\begin{center}\begin{figure}\begin{tabular}{c}
        \includegraphics[width=0.95\textwidth, height=0.16\textheight]{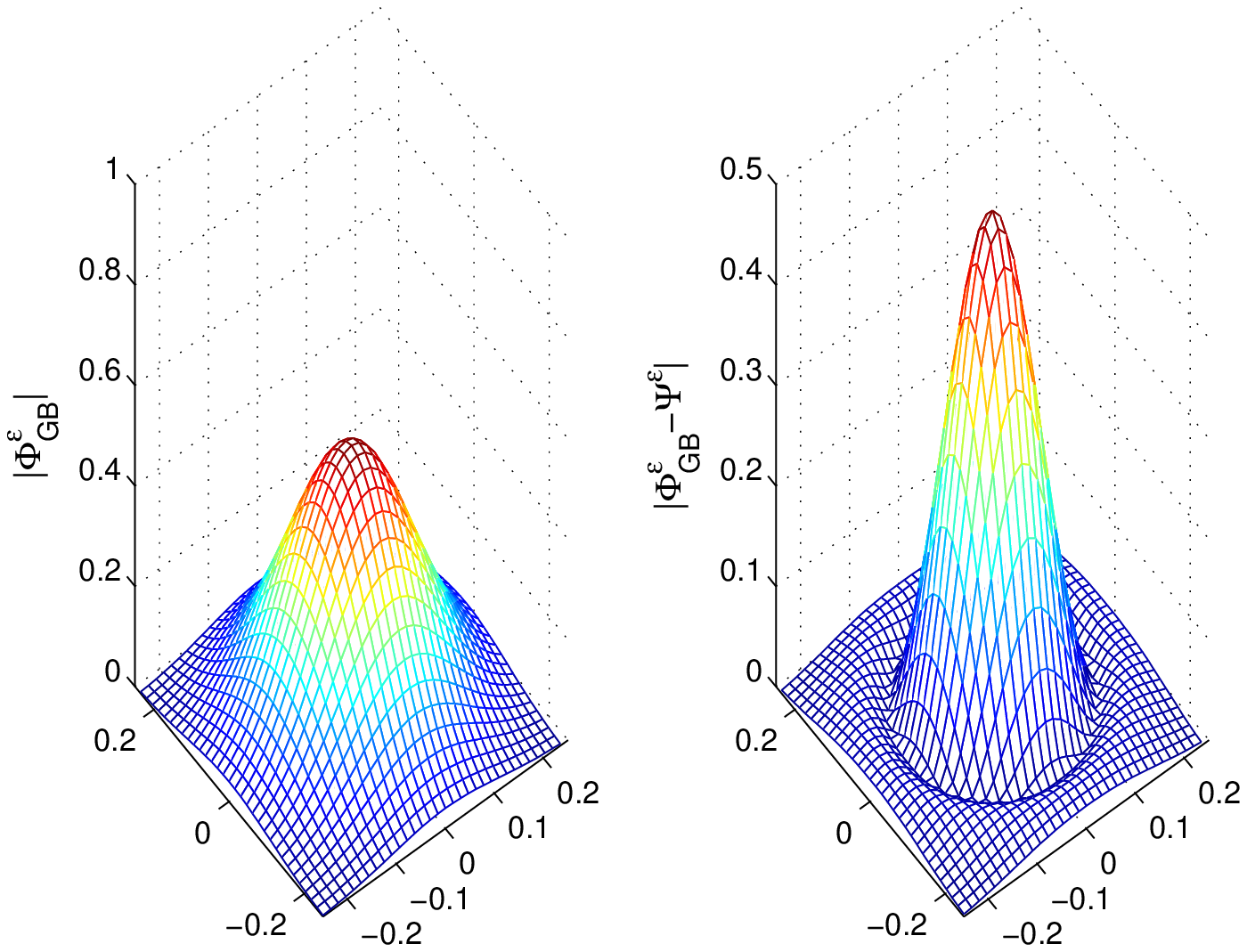} \\
        \includegraphics[width=0.95\textwidth, height=0.16\textheight]{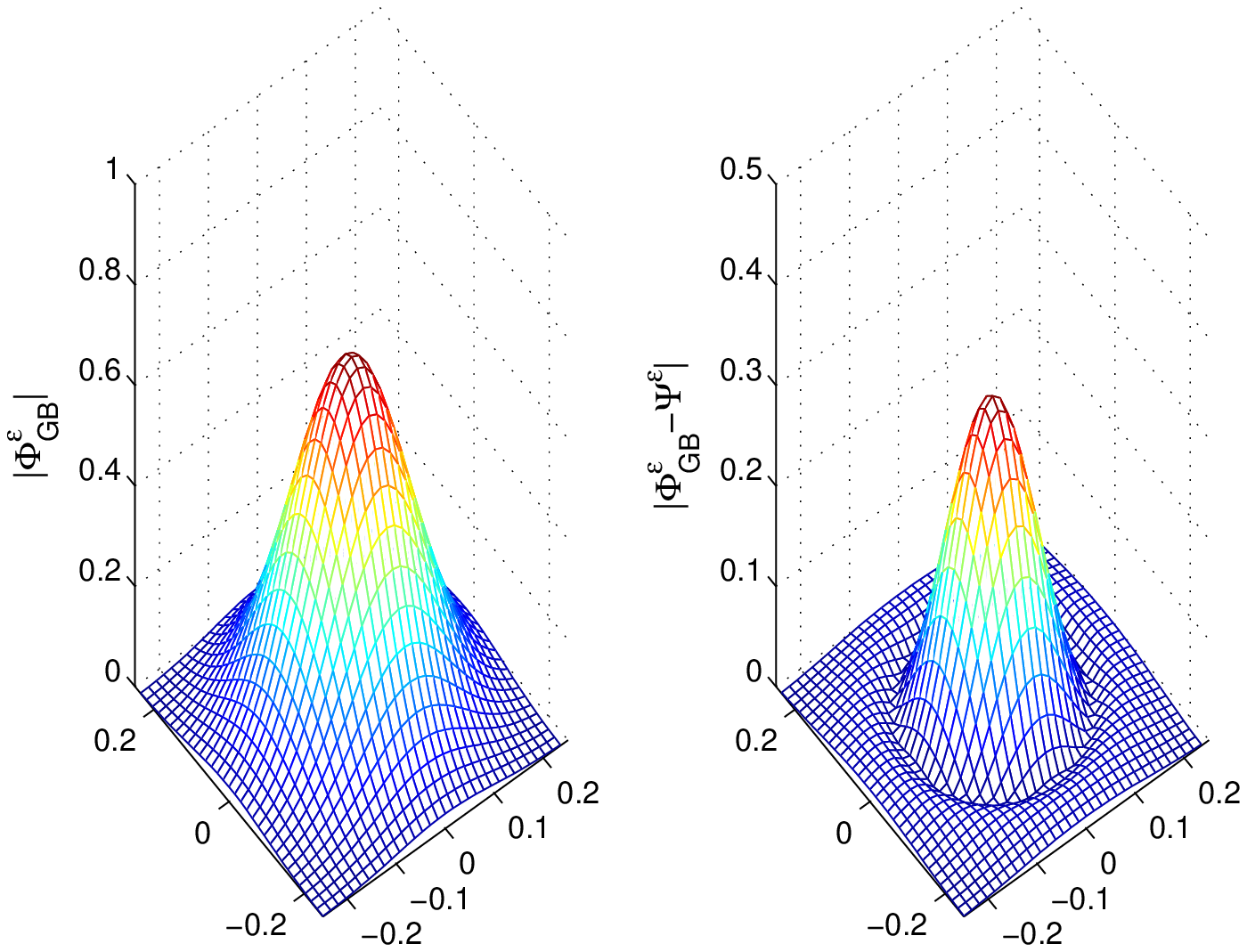} \\
        \includegraphics[width=0.95\textwidth, height=0.16\textheight]{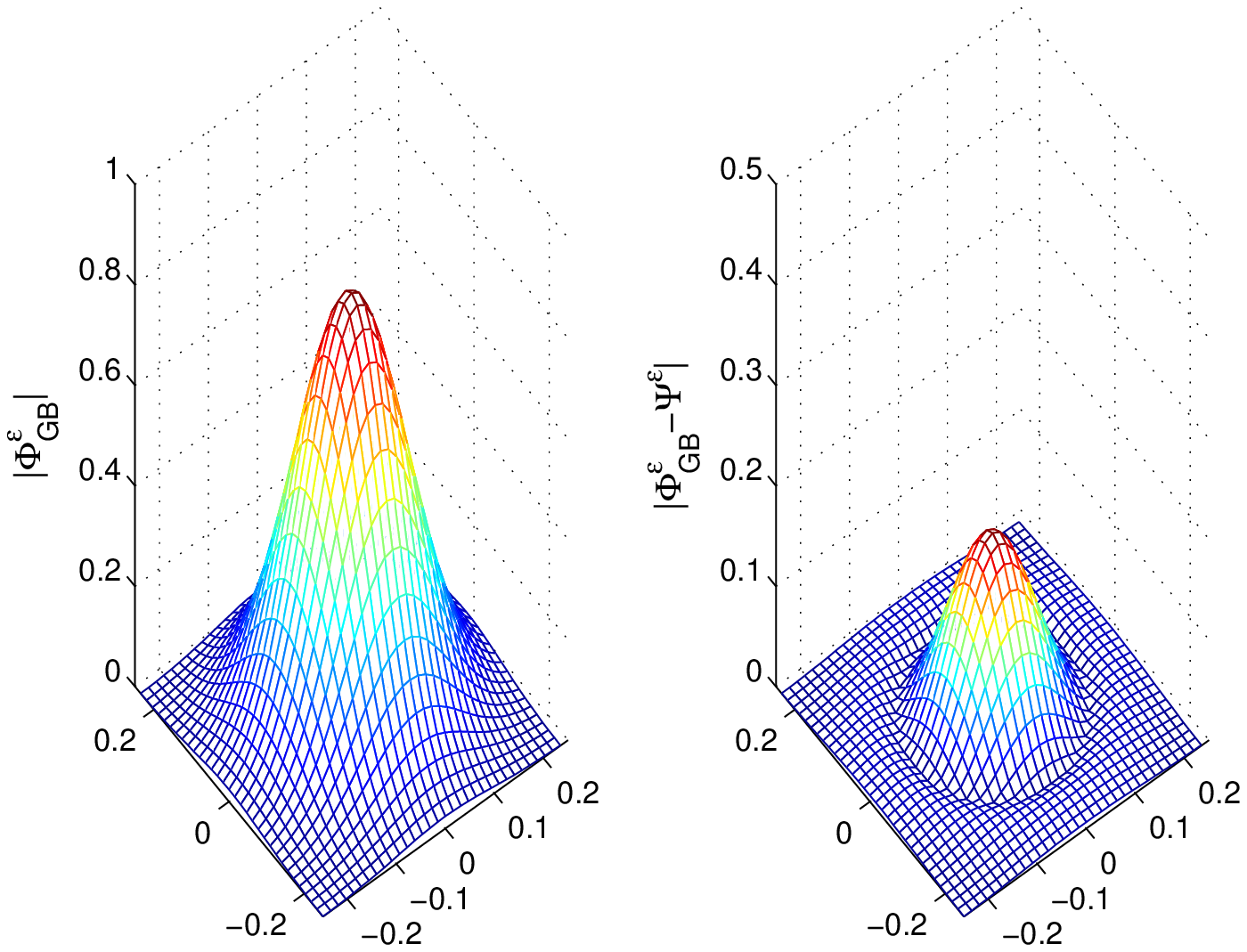}
    \end{tabular}
    \caption{Example 1, at time $t=0.5,\;x_3=0$, from top to down, they are the amplitude of Gaussian beam solutions (Left) $\abs{\bphi_{GB}^{\epsilon}}$ and the absolute error (Right) $\abs{\bpsi^{\epsilon}-\bphi_{GB}^{\epsilon}}$ for different $\epsilon=\frac{1}{256},\; \frac{1}{512},\; \frac{1}{1024}$.} \label{ex:fig11}
\end{figure}\end{center}

% ---------exam01---------
\nin {\bf Example 1.} We consider the zero external fields, i.e. $V(\bx)=0$ and $\bA(\bx)=0$. The initial condition for the Dirac equation \eqref{eqn:dirac}-\eqref{eqn:diracinit} is
\begin{equation*}
    \bpsi^{\epsilon}_0(\bx)=e^{-\frac{\abs{\bx}^2}{4d^2}}\chi, \quad \chi=(1,0,0,0)^T, \quad d=\frac{1}{16}.
\end{equation*}
In this example, $\Psi_2^{\epsilon}=\Psi_3^{\epsilon}=\Psi_4^{\epsilon}=0$ and
\begin{equation*}
    i\epsilon\partial_t \Psi_1^{\epsilon}=\Psi_1^{\epsilon},
\end{equation*}
which can be explicitly solved as
\begin{equation*}
    \Psi^{\epsilon}_1(t,\bx)=e^{-\frac{\abs{\bx}^2}{4d^2}}e^{-\frac{it}{\epsilon}}.
\end{equation*}
The $l^1,\;l^2$ and $l^{\infty}$ errors between the solutions of the Dirac equation $\bpsi^{\epsilon}$ and those of the Gaussian beam method $\bphi^{\epsilon}_{GB}$ for different $\epsilon$ are given in Table \ref{ex:table11}. Here we take $t=0.5$, the time step and mesh size of the Gaussian beam method satisfy $\Delta t=O(\epsilon^{\half}),\;\Delta y=O(\epsilon^{\half})$. We plot the wave amplitudes and absolute errors for different $\epsilon$ in Figure \ref{ex:fig11}, from which, one can see the Gaussian beam methods is more accurate for small $\epsilon$ and converges nearly first order with respect of $\epsilon$. On the other hand, the absolute error could be large for big $\epsilon$, e.g. the absolute $l^{\infty}$ error of the Gaussian beam solution could be $0.500$ for $\epsilon=\frac{1}{256}$. It is because we are considering the asymptotic numerical method. The approximations may not good when the asymptotic parameter $\epsilon$ is not small enough. However, the $l^{\infty}$ error decays almost linearly in $\epsilon$, and we can still conclude that the Gaussian beam method is accurate and efficient in the semi-classical regime $\epsilon\ll1$.

\begin{table}\begin{center}\begin{tabular}{ccccc}
        \hline $\epsilon$ & $\frac{1}{256}$ & $\frac{1}{512}$ & $\frac{1}{1024}$ & $\frac{1}{2048}$\\ \hline \\[-2mm]
        $\norm{\bpsi^{\epsilon}-\bphi^{\epsilon}_{GB}}_1$ & $7.08\times10^{-2}$ & $4.15\times10^{-2}$ & $2.27\times10^{-2}$ & $1.19\times10^{-2}$ \\
        $\norm{\bpsi^{\epsilon}-\bphi^{\epsilon}_{GB}}_2$ & $1.29\times10^{-1}$ & $7.79\times10^{-2}$ & $4.31\times10^{-2}$ & $2.30\times10^{-2}$ \\
        $\norm{\bpsi^{\epsilon}-\bphi^{\epsilon}_{GB}}_{\infty}$ & $5.00\times10^{-1}$ & $3.16\times10^{-1}$ & $1.80\times10^{-1}$ & $9.71\times10^{-2}$ \\ \hline
    \end{tabular}\end{center}
    \caption{The $l^1,\;l^2$ and $l^{\infty}$ errors of the solutions at $t=0.50$ for Example 1. The convergence rate in $\epsilon$ are $0.8591$ in the $l^1$ norm, $0.8311$ in the $l^2$ norm and $0.7912$ in the $l^{\infty}$ norm respectively.} \label{ex:table11}
\end{table}

\begin{center}\begin{figure}\begin{tabular}{c}
        \includegraphics[width=0.95\textwidth, height=0.16\textheight]{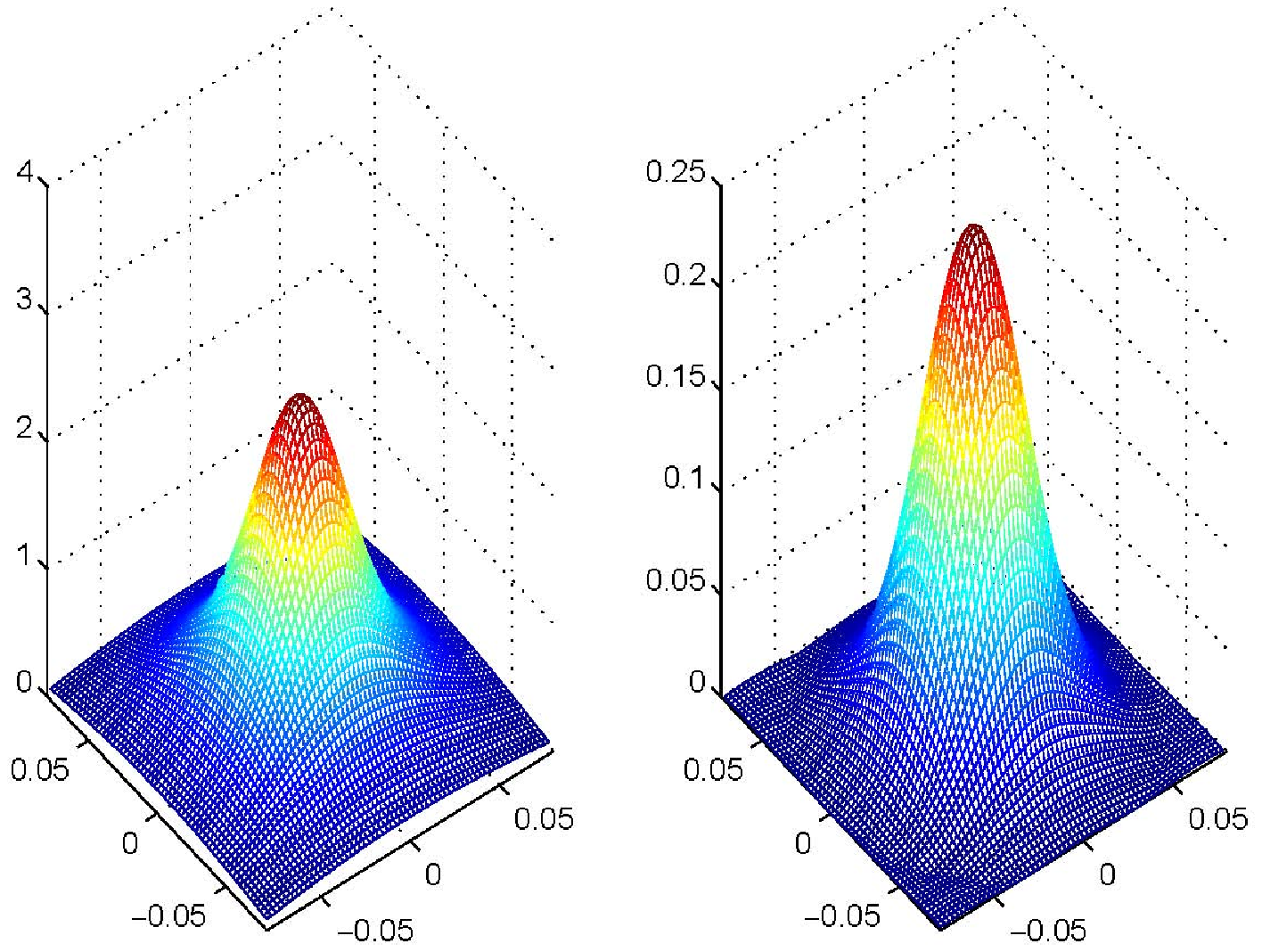} \\
        \includegraphics[width=0.95\textwidth, height=0.16\textheight]{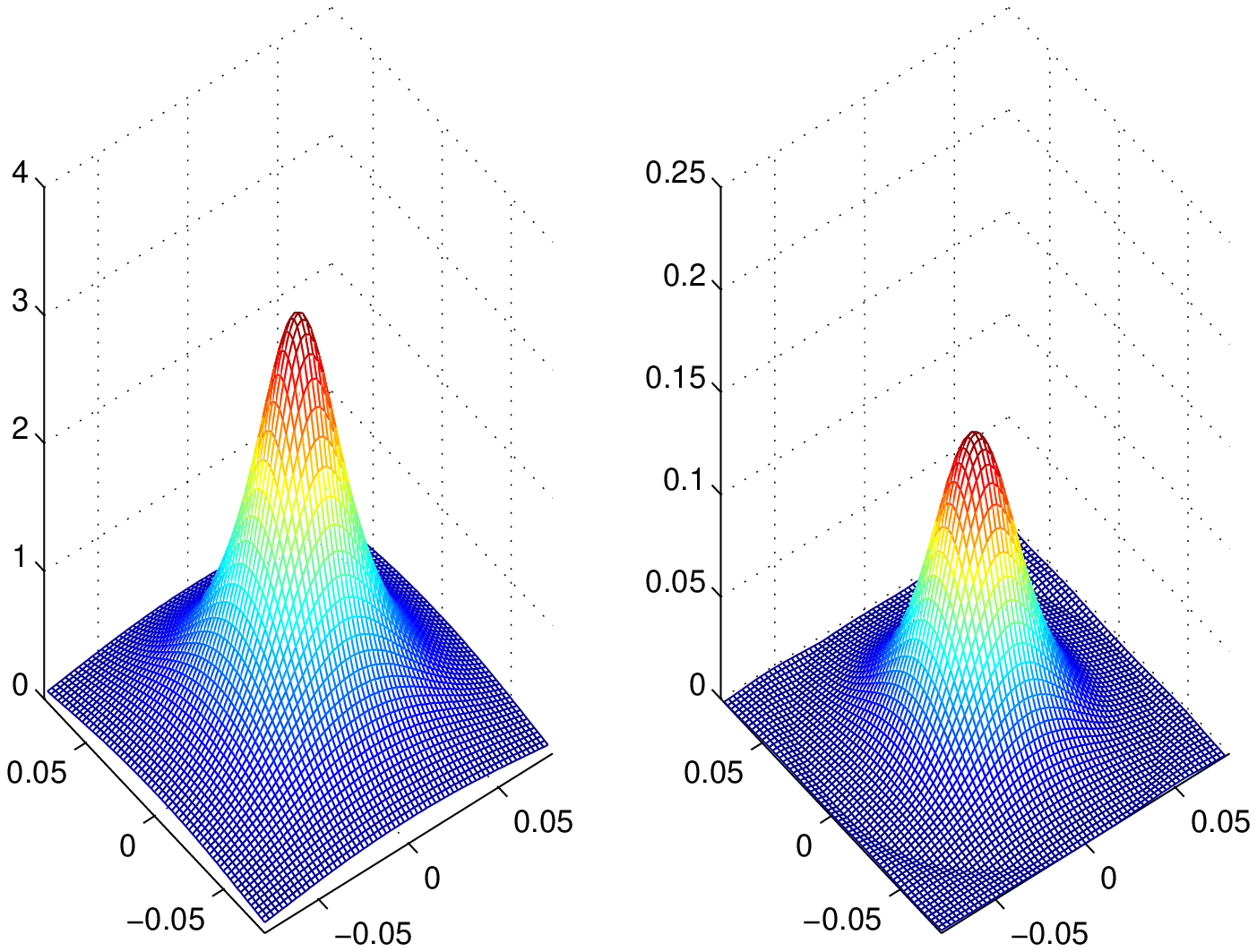} \\
        \includegraphics[width=0.95\textwidth, height=0.16\textheight]{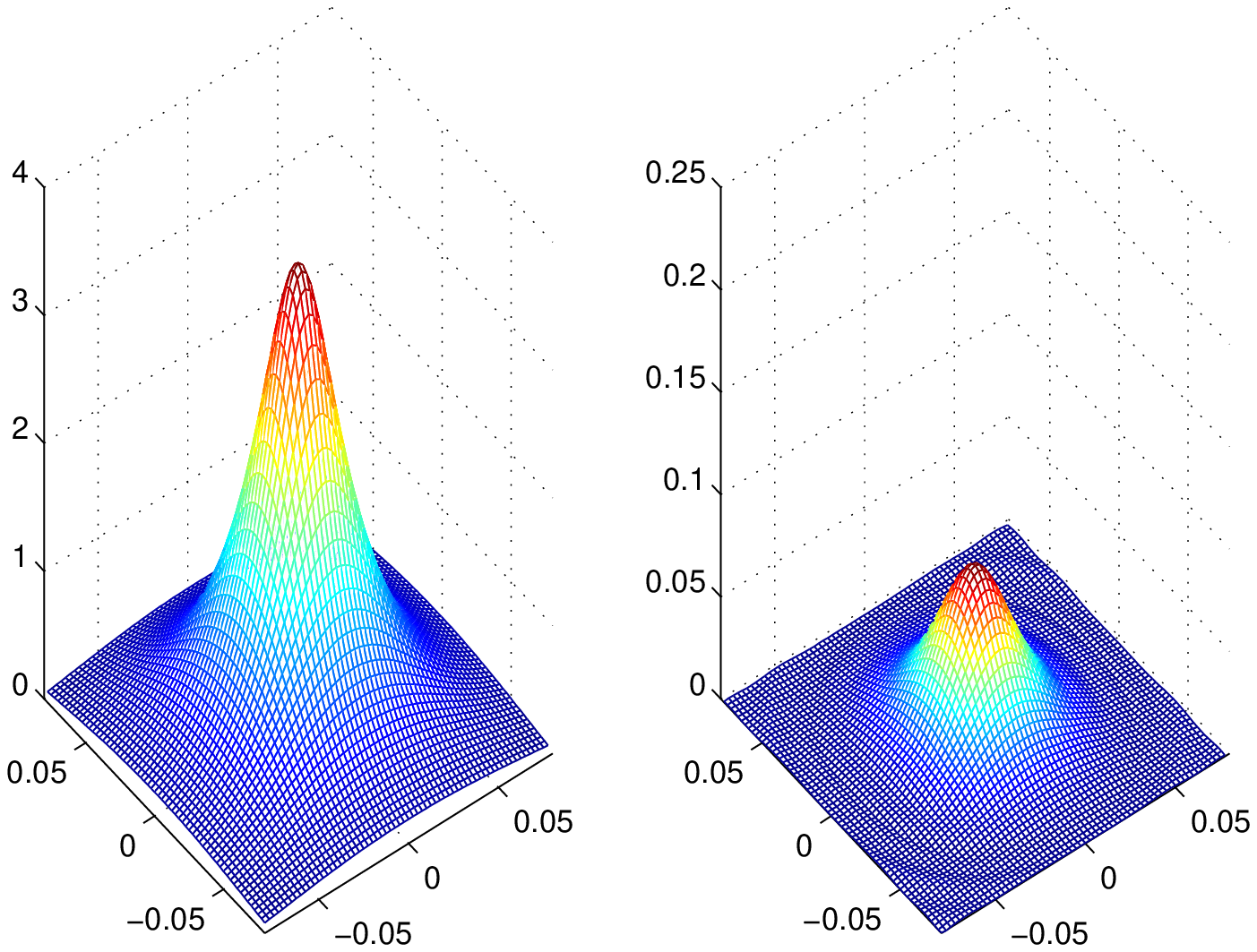}
    \end{tabular}
    \caption{Example 2, at time $t=0.38,\; x_3=0$, from top to down, they are the amplitude of Gaussian beam solutions (Left) $\abs{\bphi_{GB}^{\epsilon}}$ and the relative error (Right) $\frac{\abs{\bpsi^{\epsilon}-\bphi_{GB}^{\epsilon}}}{\norm{\bpsi^{\epsilon}}_{\infty}}$ for different $\epsilon=\frac{1}{512},\; \frac{1}{1024},\; \frac{1}{2048}$.} \label{ex:fig21}
\end{figure}\end{center}

\begin{center}\begin{figure}\begin{tabular}{c}
        \includegraphics[width=0.95\textwidth, height=0.16\textheight]{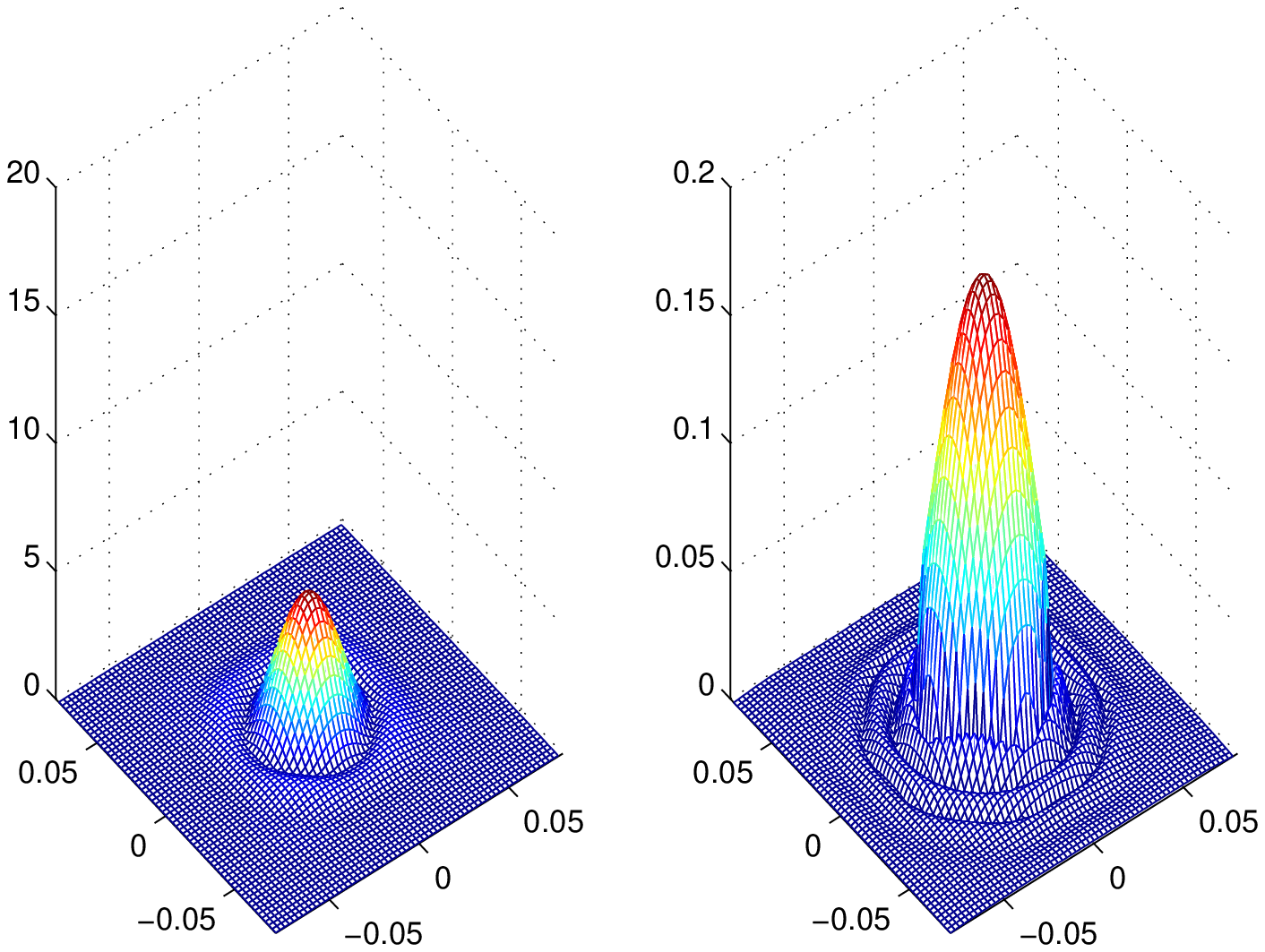} \\
        \includegraphics[width=0.95\textwidth, height=0.16\textheight]{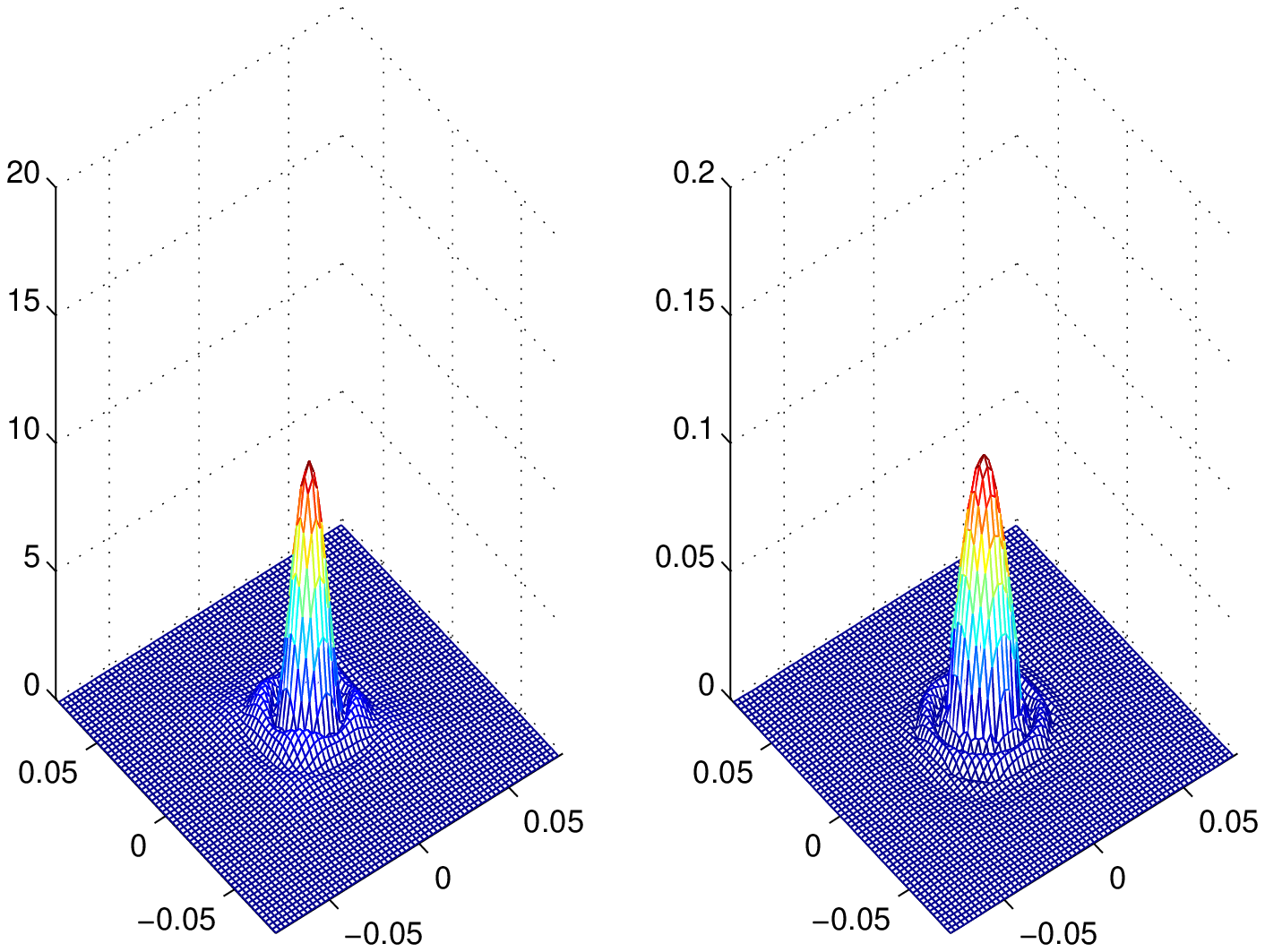} \\
        \includegraphics[width=0.95\textwidth, height=0.16\textheight]{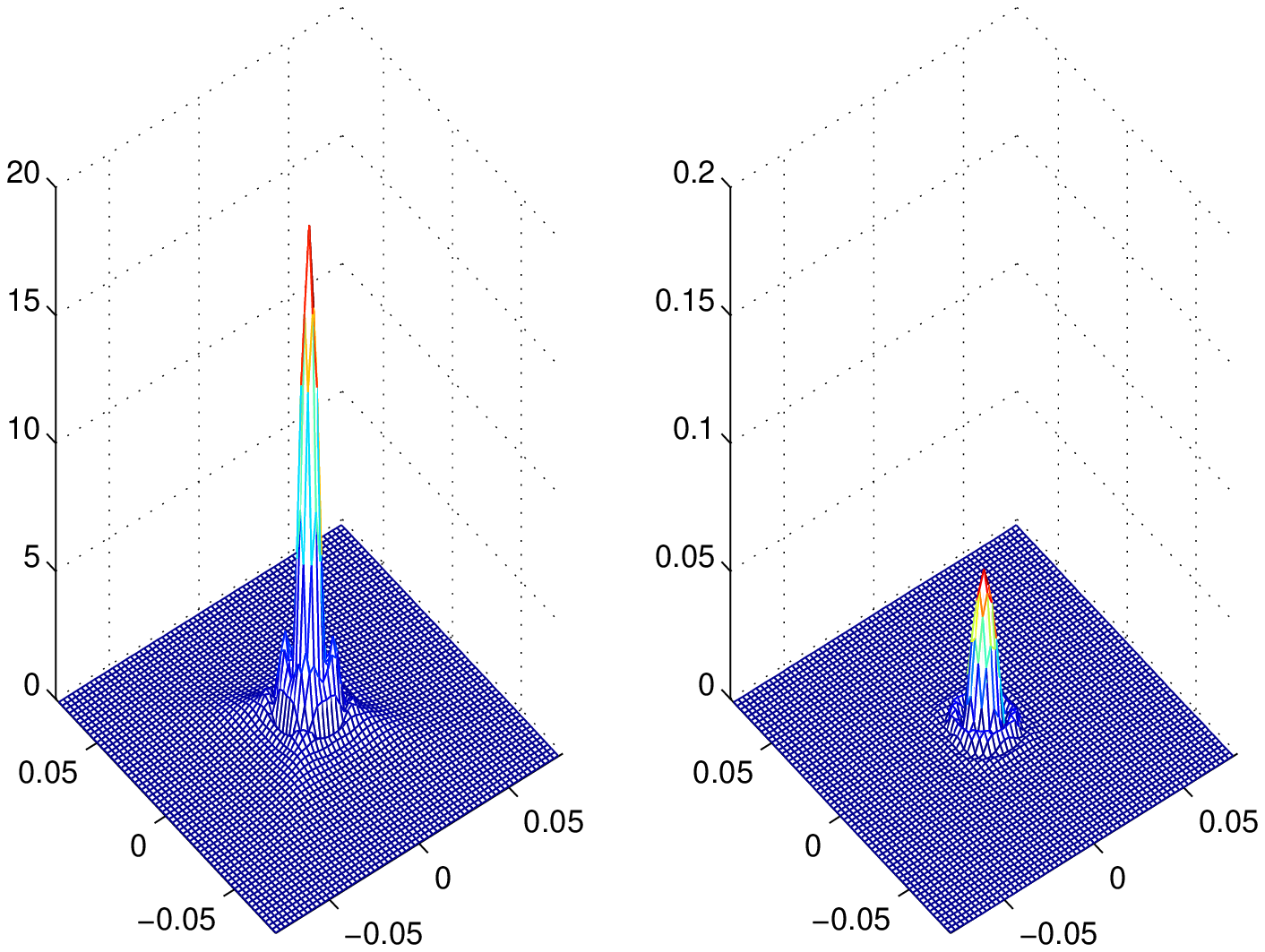}
    \end{tabular}
    \caption{Example 2, at time $t=0.56,\; x_3=0$, from top to down, they are the amplitude of Gaussian beam solutions (Left) $\abs{\bphi_{GB}^{\epsilon}}$ and the relative error (Right) $\frac{\abs{\bpsi^{\epsilon}-\bphi_{GB}^{\epsilon}}}{\norm{\bpsi^{\epsilon}}_{\infty}}$ for different $\epsilon=\frac{1}{512},\; \frac{1}{1024},\; \frac{1}{2048}$.} \label{ex:fig22}
\end{figure}\end{center}

% ---------exam02---------
\nin {\bf Example 2.} We consider the same zero external fields, the initial condition for the Dirac equation \eqref{eqn:dirac}-\eqref{eqn:diracinit} is
\begin{eqnarray*}
    \bpsi^{\epsilon}_0(\bx) &=& e^{-\frac{\abs{\bx}^2}{4d^2}}e^{iS_0(\bx)/\epsilon}\chi(\bx), \quad d=\frac{1}{16}, \\
    S_0(\bx) &=& \frac{1}{40}(1+\cos 2\pi x_1)(1+\cos 2\pi x_2), \\
    \chi(\bx) &=& \brac{\half\brac{\sqrt{(\partial_{x_1}S_0)^2+(\partial_{x_2}S_0)^2+1}+1},
        0, 0, \half\brac{\partial_{x_1}S_0+\partial_{x_2}S_0}}^T.
\end{eqnarray*}
In this example, $\Psi^{\epsilon}_2=\Psi^{\epsilon}_3=0$ and
\begin{eqnarray*}
    i\epsilon\partial_t \Psi_1 &=& \Psi_1-i\epsilon\partial_{x_1}\Psi_4-\epsilon\partial_{x_2}\Psi_4, \\
    i\epsilon\partial_t \Psi_4 &=& -i\epsilon\partial_{x_1}\Psi_1+\epsilon\partial_{x_2}\Psi_1-\Psi_4,
\end{eqnarray*}
which reduces to a two dimensional problem and can be solved by the time-splitting spectral method in only one time step. Due to the compressive initial velocity, the caustics will form at about $t\approx 0.56$. The $l^1,\;l^2$ and $l^{\infty}$ errors between the solutions of the Dirac equation $\bpsi^{\epsilon}$ and those of the Gaussian beam method $\bphi^{\epsilon}_{GB}$ for different $\epsilon$ are given in Tables \ref{ex:table21}-\ref{ex:table22}, for $t=0.375$ and $t=0.56$ respectively. We remark that in Table \ref{ex:table22} we compare the relative $l^{\infty}$-error since the caustics form. We plot the wave amplitudes and relative errors for different $\epsilon$ in Figures \ref{ex:fig21}-\ref{ex:fig22}, for which, we can draw the same conclusions as in Example 1.

\begin{table}\begin{center}\begin{tabular}{ccccc}
        \hline $\epsilon$ & $\frac{1}{512}$ & $\frac{1}{1024}$ & $\frac{1}{2048}$ & $\frac{1}{4096}$ \\ \hline \\[-2mm]
        $\norm{\bpsi^{\epsilon}-\bphi^{\epsilon}_{GB}}_1$ & $1.23\times10^{-1}$ & $6.41\times10^{-2}$ & $3.26\times10^{-2}$ & $1.71\times10^{-2}$ \\
        $\norm{\bpsi^{\epsilon}-\bphi^{\epsilon}_{GB}}_2$ & $2.21\times10^{-1}$ & $1.23\times10^{-1}$ & $6.53\times10^{-2}$ & $3.37\times10^{-2}$ \\
        $\norm{\bpsi^{\epsilon}-\bphi^{\epsilon}_{GB}}_{\infty}$ & $8.36\times10^{-1}$ & $5.40\times10^{-1}$ & $3.09\times10^{-1}$ & $1.57\times10^{-1}$ \\ \hline
    \end{tabular}\end{center}
    \caption{The $l^1,\;l^2$ and $l^{\infty}$ errors of the solutions at $t=0.38$ for Example 2. The convergence rate in $\epsilon$ are $0.9490$ in the $l^1$ norm, $0.9051$ in the $l^2$ norm and $0.8111$ in the $l^{\infty}$ norm respectively.} \label{ex:table21}
\end{table}

\begin{table}\begin{center}\begin{tabular}{ccccc}
        \hline $\epsilon$ & $\frac{1}{512}$ & $\frac{1}{1024}$ & $\frac{1}{2048}$ & $\frac{1}{4096}$ \\ \hline \\[-2mm]
        $\norm{\bpsi^{\epsilon}-\bphi^{\epsilon}_{GB}}_1$ & $8.62\times10^{-2}$ & $2.92\times10^{-2}$ & $1.14\times10^{-2}$ & $5.11\times10^{-3}$ \\
        $\norm{\bpsi^{\epsilon}-\bphi^{\epsilon}_{GB}}_2$ & $2.28\times10^{-1}$ & $1.26\times10^{-1}$ & $6.62\times10^{-2}$ & $3.49\times10^{-2}$ \\
        $\frac{\norm{\bpsi^{\epsilon}-\bphi^{\epsilon}_{GB}}_{\infty}}{\norm{\bphi^{\epsilon}_{GB}}_{\infty}}$ & $1.76\times10^{-1}$ & $1.06\times10^{-1}$ & $6.14\times10^{-2}$ & $3.23\times10^{-2}$ \\ \hline
    \end{tabular}\end{center}
    \caption{The $l^1,\;l^2$ and $l^{\infty}$ errors of the solutions at $t=0.56$ for Example 2. The convergence rate in $\epsilon$ are $1.3682$ in the $l^1$ norm, $0.9030$ in the $l^2$ norm and $0.8177$ in the $l^{\infty}$ norm (relative errors) respectively.} \label{ex:table22}
\end{table}

\begin{center}\begin{figure}
    \includegraphics[width=1\textwidth]{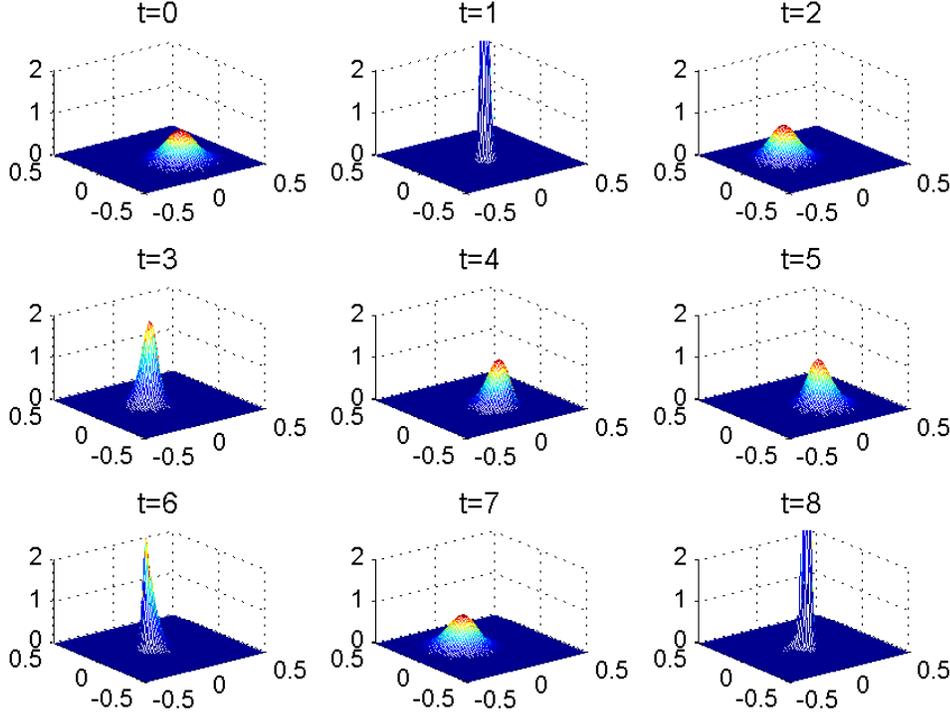}
    \caption{Example 3, the amplitude of Gaussian beam solutions at different time, here $\epsilon=\frac{1}{512},\; x_3=0$. Note at time $t=1$ and $t=8$ the amplitude are cutted since it is too big near caustics.} \label{ex:fig31}
\end{figure}\end{center}

% ---------exam03---------
\nin {\bf Example 3 (Harmonic oscillator).} We consider the zero external magnetic potential $\bA(\bx)=0$ and the quadratic external electric potential $V(\bx)=\half\abs{\bx}^2$. The initial condition for the Dirac equation \eqref{eqn:dirac}-\eqref{eqn:diracinit} is
\begin{equation*}
    \bpsi^{\epsilon}_0(\bx)=e^{-\frac{(x_1-0.1)^2+(x_2+0.1)^2+x_3^2}{4d^2}}\chi, \quad \chi=(1,0,0,0)^T, \quad d=\frac{1}{16}.
\end{equation*}
This is a full 3D problem. The time splitting spectral method is very expensive because of the large requirement of memory for very small $\epsilon$, which the Gaussian beam method is accurate. In Figure \ref{ex:fig31}, we depict the wave amplitude at different time $t$. In this example, we choose $\epsilon=\frac{1}{512}$. We can see that the wave packet moves in circles due to its interaction with the harmonic external potential.

\medskip

%*****************************************************************************************
%************* Conclusion ****************************************************************
%*****************************************************************************************

\section{Conclusion} \label{sec:con}\setcounter{equation}{0}
In this work, we developed the Gaussian beam method for the Dirac equation. The Eulerian Gaussian beam method provides a simple way to compute the Hessian matrices for the phase, as in \cite{JiWuYa:08}. The proposed method is shown numerically to be accurate and efficient. The required mesh size and time step should be of $O(\sqrt{\epsilon})$. Compare to the traditional numerical method, which requires mesh size to be $O(\epsilon)$, the computational cost for our method is cheap, especially when $\epsilon$ is very small. A more interesting question is to simulate the graphene by using the Gaussian beam method. We are currently investigating this important model and hope to report our progress in the near future.

\section*{Appendix}
\renewcommand{\theequation}{A-\arabic{equation}}

In this appendix, we give the detailed derivation of the Lagrangian Gaussian beam method for the Dirac equation. For convenience, we drop the superscript $\pm$ in \eqref{eqn:LGB}
\begin{equation}
    \boldsymbol{\phi}^{\epsilon}_{la}(t,\bx,\by_0)=\bu(t,\bx,\by)e^{iT(t,\bx,\by)/\epsilon},
\end{equation}
with
\begin{eqnarray}
    T(t,\bx,\by) &=& S(t,\by)+\bxi(t,\by)\cdot(\bx-\by)+\half(\bx-\by)^TM(t,\by)(\bx-\by), \label{eqn:gb_Texpan} \\
    \bu(t,\bx,\by) &=& \widetilde{\bu}_0(t,\bx,\by)+\epsilon\widetilde{\bu}_1(t,\bx,\by), \label{eqn:gb_uexpan} \\
    \widetilde{\bu}_0(t,\bx,\by) &=& \bu_0(t,\by)+(\bx-\by)\cdot\nabla_{\by}\bu_0(t,\by)
        +\half(\bx-\by)^T\brac{\nabla_{\by\by}\bu_0(t,\by)}(\bx-\by), \label{eqn:gb_u0expan} \\
    \widetilde{\bu}_1(t,\bx,\by) &=& \bu_1(t,\by). \nn
\end{eqnarray}
%Then we have
%\begin{eqnarray*}
%    \partial_t\phi &=& \brac{\partial_t(\bu_0+\epsilon\bu_1)
%        +\frac{\rd \by}{\rd t}\cdot\nabla_{\by}(\bu_0+\epsilon\bu_1)}\exp\{\cdot\} \\
%        && +\frac{i}{\epsilon}\brac{\partial_t T+\frac{\rd \by}{\rd t}
%        \cdot\nabla_{\by}T}(\bu_0+\epsilon\bu_1)\exp\{\cdot\}, \\
%    \balpha\cdot\nabla_{\bx}\phi &=& \balpha\cdot\nabla_{\bx}(\bu_0+\epsilon\bu_1)\exp\{\cdot\}+
%        \frac{i}{\epsilon}\balpha\cdot\nabla_{\bx}T(\bu_0+\epsilon\bu_1)\exp\{\cdot\},
%\end{eqnarray*}
%then we have
%\begin{eqnarray*}
%    0 &=& i\epsilon\brac{\partial_t \bu_0+\frac{\rd \by}{\rd t}\cdot\nabla_{\by}\bu_0}
%        -\brac{\partial_t T+\frac{\rd \by}{\rd t}
%        \cdot\nabla_{\by}T}(\bu_0+\epsilon\bu_1) \\
%    && +i\epsilon \balpha\cdot\nabla_{\bx}\bu_0
%        -\balpha\cdot\nabla_{\bx}T(\bu_0+\epsilon\bu_1)
%        -\brac{-\balpha\cdot\bA+\beta+V}\brac{\bu_0+\epsilon\bu_1}.
%\end{eqnarray*}
Without loss of generality, we can assume that $\Pi(\bx,\nabla_{\by} S)\bu_0=\bu_0$.
Taking this into \eqref{eqn:dirac} and matching the first two leading order of $\epsilon$, one obtains
\begin{align}
    & \brac{\partial_tT+\frac{\rd \by}{\rd t}\cdot\nabla_{\by}T+
        \mcd(\bx,\nabla_{\bx}T)}\widetilde{\bu}_0 = 0, \label{eqn:gb_asy0} \\
    & i\brac{\partial_t\widetilde{\bu}_0+\frac{\rd \by}{\rd t}\cdot\nabla_{\by}\widetilde{\bu}_0+
    \balpha\cdot\nabla_{\bx}\widetilde{\bu}_0}
    -\brac{\partial_tT+\frac{\rd \by}{\rd t}\cdot\nabla_{\by}T+
        \mcd(\bx,\nabla_{\bx}T)}\widetilde{\bu}_1 = 0. \label{eqn:gb_asy1}
\end{align}
In order to get a nontrivial solution $\widetilde{\bu}_0\ne0$ in \eqref{eqn:gb_asy0}, we need
\begin{equation*}
    \det\brac{\partial_tT+\frac{\rd \by}{\rd t}\cdot\nabla_{\by}T+
        \mcd(\bx,\nabla_{\bx}T)}=0.
\end{equation*}
Since $h(\bx,\bxi)$ is the eigenvalue of $\mcd(\bx,\bxi)$, one gets
\begin{equation} \label{eqn:gb_S}
    \partial_tT+\frac{\rd \by}{\rd t}\cdot\nabla_{\by}T+h(\bx,\nabla_{\bx}T)=0.
\end{equation}
Taking the first and second order derivatives with respect to $\bx$ in \eqref{eqn:gb_S} gives
\begin{align}
    & \partial_t(\nabla_{\bx}T)+\frac{\rd \by}{\rd t}\cdot\nabla_{\by\bx}T
        +\nabla_{\bx}h+\nabla_{\bxi}h\cdot\nabla_{\bx\bx}T = 0, \label{eqn:gb_xi} \\
    & \partial_t(\nabla_{\bx\bx}T)+\frac{\rd \by}{\rd t}\cdot\nabla_{\by\bx\bx}T
        +\nabla_{\bx\bx}h+\nabla_{\bx\bxi}h \nabla_{\bx\bx}T \nn \\
    & \qquad +\nabla_{\bx\bx}T \nabla_{\bxi\bx}h+\nabla_{\bx\bx}T \nabla_{\bxi\bxi}h \nabla_{\bx\bx}T
        +\nabla_{\bxi}h\cdot\nabla_{\bx\bx\bx}T = 0. \label{eqn:gb_M}
\end{align}
Consider \eqref{eqn:gb_Texpan} and evaluating \eqref{eqn:gb_asy1}-\eqref{eqn:gb_M} at $\bx=\by$ yield
\begin{align}
    & \partial_t S+\frac{\rd \by}{\rd t}\cdot\brac{\nabla_{\by}S-\bxi}+h(\bx,\bxi) = 0, \label{eqn:gb_S2} \\
    & \partial_t \bxi+\frac{\rd \by}{\rd t}\cdot\brac{\nabla_{\by}\bxi-M}+\nabla_{\by}h+
        \nabla_{\bxi}h\cdot M=0,  \label{eqn:gb_xi2} \\
    & \partial_t M+\frac{\rd \by}{\rd t}\cdot\nabla_{\by}M+\nabla_{\by\by}h
        +\nabla_{\by\bxi}hM+M\nabla_{\bxi\by}h+M\nabla_{\bxi\bxi}hM=0, \label{eqn:gb_M2} \\
    & i\brac{\partial_t\bu_0+\balpha\cdot\nabla_{\by}\bu_0}+
    \brac{\partial_t S+\frac{\rd \by}{\rd t}\cdot\brac{\nabla_{\by}S-\bxi}+\mcd(\bx,\bxi)}\bu_1 = 0. \label{eqn:gb_u}
\end{align}
We choose the beam center that satisfies
\begin{equation*}
    \frac{\rd \by}{\rd t}=\nabla_{\bxi} h,
\end{equation*}
then \eqref{eqn:gb_S2}-\eqref{eqn:gb_M2} can be written as
\begin{eqnarray*}
    \frac{\rd \bxi}{\rd t} &=& -\nabla_{\by}h, \\
    \frac{\rd  S}{\rd t} &=& \nabla_{\bxi}h\cdot\bxi-h, \\
    \frac{\rd M}{\rd t} &=& -\nabla_{\by\by}h-\nabla_{\by\bxi}hM
        -M\nabla_{\bxi\by}h-M\nabla_{\bxi\bxi}hM.
\end{eqnarray*}
To obtain the transport equation for $\bu_0$, we apply the projection $\Pi(\bx,\nabla_{\by} S)$ to \eqref{eqn:gb_u}:
\begin{equation*}
    \Pi(\bx,\nabla_{\by} S)\brac{\partial_t+\balpha\cdot\nabla_{\by}}\bu_0=0.
\end{equation*}
This is the solvability condition for $\bu_1$, and one can finally get the ODE for Lagrangian Gaussian beam amplitude after a series of calculations
\begin{equation*}
    \frac{\rd \bu_0}{\rd t} = -\half \brac{\nabla_{\by}\cdot\omega}\bu_0+
        \mca \bu_0.
\end{equation*}

\medskip
\begin{rem}
    The amplitude \eqref{eqn:gb_uexpan}-\eqref{eqn:gb_u0expan} is expanded to higher order in both the Taylor expansion and the asymptotic expansion, which is different from the Gaussian beam method for the Schr\"odinger equation \cite{JiWuYa:08}. The reason is that a higher order asymptotic expansion is needed for deriving the transport equation for the amplitude when using the solvability condition. The higher order Taylor expansion should be used to match the high order asymptotic expansion for the Gaussian beam method.
\end{rem}

\medskip
\begin{rem}
    The equation \eqref{eqn:gb_u0expan} can be written in a more general form, e.g.
    \begin{equation*}
        \widetilde{\bu}_0(t,\bx,\by) = \bu_{00}(t,\by)+(\bx-\by)\cdot \bu_{01}(t,\by)
            +\half(\bx-\by)^T \bu_{02}(t,\by)(\bx-\by).
    \end{equation*}
    where $\bu_{00}\in C^{\infty}(\br^4,\bc^4),\;\bu_{01}\in C^{\infty}(\br^4,\bc^{3\times4})$ and
    $\bu_{02}\in C^{\infty}(\br^4,\bc^{(3\times 3)\times4})$. Since there are more freedoms than restrictions, one can easily formulate them as
    \begin{equation*}
        \bu_{01}=\nabla_{\by}\bu_{00}, \quad
        \bu_{02}=\nabla_{\by\by}\bu_{00},
    \end{equation*}
    to close the system. This is consistent to the equation \eqref{eqn:gb_u0expan}.
\end{rem}

%{\bf Acknowledgement.}
%\medskip


\begin{thebibliography}{10}

\bibitem{CsDa:06} J. Cserti and G. D\'avid, {\em Unified description of Zitterbewegung for spintronic, graphene, and superconducting systems}, Phys. Rev. B, 74, 172305, 2006.

\bibitem{Di:28} P. Dirac, {\em The Quantum Theory of the Electron}, Proc. R. Soc. Lond. A, 117, 610-624, 1928.

\bibitem{EnRu:03} B. Engquist and O. Runborg, {\em Computational high frequency wave propagation}, Acta Numer., 12, 181-266, 2003.

\bibitem{He:75} E.J. Heller, {\em Time-dependent approach to semiclassical dynamics}, J. Chem. Phys., 62, 1544-1555, 1975.

\bibitem{HsQiWr:08} D. Hsieh, D. Qian, L. Wray, Y. Xia, Y.S. Hor, R.J. Cava and M.Z. Hasan, {\em A topological Dirac insulator in a quantum spin Hall phase}, Nature, 452, 970-974, 2008.

\bibitem{HuJi:05} Z.Y. Huang, S. Jin, P.A. Markowich, C. Sparber and C.X. Zheng, {\em A time-splitting spectral scheme for the Maxwell-Dirac system}, J. Comput. Phys., 208, 761-789, 2005.

\bibitem{JiLiOsTs:05} S. Jin, H. Liu, S. Osher and R. Tsai, {\em Computing multi-valued physical observables the semiclassical limit of the Schr\"odinger equations},  J. Comput. Phys., 205, 222-241, 2005.

\bibitem{JiMaSp:11} S. Jin, P. Markowich and C. Sparber, {\em Mathematical and computational methods for semiclassical Schr\"odinger equations}, Acta Numer., 20, 121-209, 2011.

\bibitem{JiOs:03} S. Jin and S. Osher, {\em A level set method for the computation of multivalued solutions to quasi-linear hyperbolic PDEs and Hamilton-Jacobi equations}, Commun. Math. Sci., 1, 575-591, 2003.

\bibitem{JiWuYa:08} S. Jin, H. Wu and X. Yang, {\em Gaussian beam methods for the Schr\"odinger equation in the semi-classical regime: Lagrangian and Eulerian formulations}, Commun. Math. Sci., 6, 995-1020, 2008.

\bibitem{JiWuYa:10} S. Jin, H. Wu and X. Yang, {\em A numerical study of the Gaussian beam methods for one-dimensional Schr\"odinger-Poisson equations}, J. Comput. Math., 28, 261-272, 2010.

\bibitem{JiWuYa:11} S. Jin, H. Wu and X. Yang, {\em Semi-Eulerian and high order Gaussian Beam methods for the Schr\"odinger equation in the Semiclassical regime}, Commun. Comput. Phys., 9, 668-687, 2011.

\bibitem{JiWuYaHu:10} S. Jin, H. Wu, X. Yang and Z. Huang, {\em Bloch Decomposition-based Gaussian Beam Method for the Schr\"odinger equation with periodic potentials},  J. Comput. Phys., 229, 4869-4883, 2010.

\bibitem{JoBlSa:88} W. Johnson, S. Blundell and J. Sapirstein, {\em Finite basis sets for the Dirac equation constructed from B splines}, Phys. Rev. A, 37, 307-315, 1988.

\bibitem{Kl:29} O. Klein, {\em Die Reflexion von Elektronen an einem Potentialsprung nach der relativistischen Dynamik von Dirac}, Z. Phys. A, 53, 157-165, 1929.

\bibitem{LeQi:09} S. Leung and J. Qian, {\em Eulerian Gaussian beams for Schr\"odinger equations in the semi-classical regime}, J. Comput. Phys., 228, 2951-2977, 2009.

\bibitem{LeQiBu:07} S. Leung, J. Qian and R. Burridge, {\em Eulerian Gaussian beams for high-frequency wave propagation}, Geophysics, 72, 61-76, 2007.

\bibitem{LiRuTa:pre} H. Liu, O. Runborg and N.M. Tanushev, {\em Error estimates for Gaussian Beam superpositions}, Math. Comp., to appear.

\bibitem{LuYa:11} J. Lu and X. Yang, {\em Frozen Gaussian approximation for high frequency wave propagation}, Commun. Math. Sci., 9, 663-683, 2011.

\bibitem{LuYa:CPAM} J. Lu and X. Yang, {\em Convergence of frozen Gaussian approximation for high frequency wave propagation}, Comm. Pure Appl. Math., to appear.

\bibitem{Mi:03} C. Min, {\em Simplicial isosurfacing in arbitrary dimension and codimension}, J. Comput. Phys., 190, 295-310, 2003.

\bibitem{MoSc:11} O. Morandi and F. Sch\"urrer, {\em Wigner model for quantum transport in graphene}, J. Phys. A, 44, 265301, 2011.

\bibitem{MoRu:10} M. Motamed and O. Runborg, {\em Taylor expansion errors in Gaussian beam summation}, Wave motion, 47, 421-439, 2010.

\bibitem{NeGu:09} A.H. Castro Neto, F. Guinea, N.M.R. Peres, K.S. Novoselov and A.K. Geim, {\em The electronic properties of graphene}, Rev. Modern Phys., 81, 109-162, 2009.

\bibitem{PeMeOs:99} D. Peng, B. Merriman, S. Osher, H. Zhao and M. Kang, {\em A PDE based fast local level set method}, J. Comput. Phys., 155, 410-438, 1999.

\bibitem{Po:82} M.M. Popov, {\em A new method of computation of wave fields using Gaussian beams}, Wave Motion, 4, 85-97, 1982.

\bibitem{QiYi:10} J. Qian and L. Ying, {\em Fast Gaussian wavepacket transforms and Gaussian beams for the Schr\"oinger equation}, J. Comput. Phys., 229, 7848-7873, 2010.

\bibitem{Ra:82} J. Ralston, {\em Gaussian beams and the propagation of singularities}, Studies in PDEs, MAA stud. Math., 23, 206-248, 1982.

\bibitem{Sc:30} E. Schr\"odinger, {\em \"Uber die kr\"aftefreie Bewegung in der relativistischen Quantenmechanik}, Sitzungsber. Preuss. Akad. Wiss., Phys. Math. Kl., 24, 418-428, 1930.

\bibitem{SpMa03:44} C. Sparber and P.A. Markowich, {\em Semiclassical asymptotics for the Maxwell-Dirac system}, J. Math. Phys., 44, 4555-4572, 2003.

\bibitem{SpMa03:45} C. Sparber and P.A. Markowich, {\em Erratum: Semiclassical asymptotics for the Maxwell-Dirac system}, J. Math. Phys., 45, 5101, 2003.

\bibitem{Sp:00} H. Spohn, {\em Semiclassical limit of the Dirac equation and spin precession}, Ann. Physics, 282, 420-431, 2000.

\bibitem{TaEnTs:09} N.M. Tanushev, B. Engquist and R. Tsai, {\em Gaussian beam decomposition of high frequency wave fields}, J. Comput. Phys., 228, 8856-8871, 2009.

\bibitem{Ta:08} N.M. Tanushev, {\em Superpositions and higher order Gaussian beams}, Commun. Math. Sci., 6, 449-475, 2008.

\bibitem{Th:92} B. Thaller, {\em The Dirac Equation}, Springer, 1992.

\bibitem{Wen:09} X. Wen, {\em High order numerical methods to two dimensional delta function integrals in level set methods}, J. Comput. Phys., 228, 4273-4290, 2009.

\bibitem{Wen:10} X. Wen, {\em High order numerical methods to three dimensional delta function integrals in level set methods}, SIAM J. Sci. Comput., 32, 1288-1309, 2010.

\bibitem{XiQi:09} Y. Xia, D. Qian, D. Hsieh, L.Wray, A. Pal, H. Lin, A. Bansil, D. Grauer, Y.S. Hor, R.J. Cava and M.Z. Hasan, {\em Observation of a large-gap topological-insulator class with a single Dirac cone on the surface}, Nature Physics, 5, 398-402, 2009.

\bibitem{YiZh:11} D. Yin and C. Zheng, {\em Gaussian beam formulations and interface conditions for the one-dimensional linear Schr\"odinger equation}, Wave Motion, 48, 310-324, 2011.

\bibitem{YiZh:12} D. Yin and C. Zheng, {\em Composite coherent states approximation for one-dimensional multi-phased wave functions}, Commun. Comput. Phys., 11, 951-984, 2012.

\end{thebibliography}
\end{document}